\documentclass[12pt,reqno]{amsart} 
\usepackage{amscd,amssymb}

\theoremstyle{plain}
\newtheorem{theorem}[subsection]{Theorem}
\newtheorem{lemma}[subsection]{Lemma}
\newtheorem{prop}[subsection]{Proposition}
\newtheorem{cor}[subsection]{Corollary}

\theoremstyle{definition}
\newtheorem{remark}[subsection]{Remark}
\newtheorem{definition}[subsection]{Definition}
\newtheorem{example}[subsection]{Example}
\newtheorem*{ack}{Acknowledgment}
\newtheorem{question}[subsection]{Question}

\newenvironment{alphenum}{

\begin{enumerate}}{\end{enumerate}}

\def\id {{\cdot\operatorname{Id}}}

\def\Q {{\mathbb Q}}
\def\C {{\mathbb C}}
\def\Z {{\mathbb Z}}
\def\ZZ {{{\mathbb Z}_2}}
\def\R {{\mathbb R}}
\def\l {p}

\def\vol{{\it vol}}
\def\length{{\it length}}

\def\sys{{\it sys}}
\newcommand{\ie}{{\it i.e.\ }}
\newcommand{\cf}{{\it cf.\ }}
\newcommand{\eg}{{\it e.g.\ }}
\newcommand{\mero}{{ \, -\! \to\, }}
\def\im{\operatorname{im}}

\begin{document}

\title[Systolic freedom of loop space] %
{Systolic freedom of loop space}

\author[M.~Katz]{Mikhail G.~Katz$^1$} \address{Department of Mathematics
and Computer Science, Bar Ilan University, Ramat Gan, 52900, Israel.
D\'epartement de Math\'ematiques,
Universit\'e Henri Poincar\'e, Vandoeuvre
54506, France} \email{katzmik@macs.biu.ac.il}
\urladdr{http://www.cs.biu.ac.il/\~{}katzmik}

\author[A.~Suciu]{Alexander I.~Suciu}
\address{Department of Mathematics,
Northeastern University,
Boston, MA 02115, USA}
\email{alexsuciu@neu.edu}
\urladdr{http://www.math.neu.edu/\~{}suciu}

%starting from version sent to Gafa on 3 April 2000 
%updated on 18 June 2001

\thanks{This is an updated version of a paper that appeared in 
{\it Geometric and Functional Analysis} {\textbf{11}} (2001), 60--73.
Visible at {\ttfamily
http://www.cs.biu.ac.il/\~{}katzmik/publications.html}}

\thanks{$^1$Supported by The Israel Science Foundation 
(grant no.~620/00--10.0).  Partially supported by the 
Emmy Noether Research Institute and the 
Minerva Foundation of Germany.}


\begin{abstract}
Given a pair of integers $m$ and $n$ such that $1 < m < n$, we show
that every $n$-dimensional manifold admits metrics of arbitrarily
small total volume, and possessing the following property: every
$m$-dimensional submanifold of less than unit $m$-volume is
necessarily torsion in homology.  

This result is different from the
case of a pair of complementary dimensions, for which a direct
geometric construction works and gives the analogous theorem in
complete generality.  In the present paper, we use Sullivan's 
telescope model for the rationalisation of a space to observe 
systolic freedom.
\end{abstract}

\maketitle

\section{Introduction}
\label{sec:intro}

Does the total volume of a Riemannian manifold impose a constraint
upon the volume of its submanifolds?  This question has interested
differential geometers since the work of M. Berger (see \cite{Be1} and
\cite[p.~192]{Be2}) in the seventies, and of M. Gromov \cite{G1},
\cite{G2} in the eighties.  In this note we answer this question in
the negative (for all but 1-dimensional submanifolds), as follows.

\begin{theorem} 
\label{thm:m}
Let $m$ and $n$ be integers such that $2\leq m<n$.  Let $X$ be an
$n$-dimensional compact smooth manifold.  Then $X$ admits metrics of
arbitrarily small total volume, such that every $m$-dimensional
orientable submanifold of less than unit $m$-volume is null-homologous
as a cycle with rational coefficients.
\end{theorem}

In other words, the manifold $X$ admits a sequence of metrics
$(X,g_j)$ with $\vol(g_j)\to 0$ as $j\to \infty$, and such that for
every orientable $M\subset X$ we have
\begin{equation*}
\vol_m(M)<1\implies [M]=0 \text{ in } H_m(X,\Q),
\end{equation*}
where $\vol_m(M)$ is the volume induced by $g_j$.  Such a phenomenon
is called {\em systolic freedom (modulo torsion)}.  For $m=1$ the
theorem is false: the length of the shortest noncontractible loop is
typically constrained by total volume (see section \ref{sec:i}).

Our theorem generalizes a previous result by I. Babenko and the
authors \cite{BKS}, concerning {systolic freedom} in middle dimension
(the case $n=2m$).  The next step toward the proof of
Theorem~\ref{thm:m} was accomplished in \cite{KS}, where we reduced
the problem to the examination of a finite list of CW complexes (the
notion of systolic freedom can be extended in this setting).  These
complexes are the successive skeleta (up to the $n^{th}$ skeleton) of
the based loop space, $\Omega(S^{m+1})$, of the sphere $S^{m+1}$.

The main idea of the present paper is to use rational homotopy theory
to establish the systolic freedom of the skeleta of the loop space.

The starting point of our proof is a map from the $p$-fold Cartesian
product $(S^m)^{\times p}$, where $p=[\frac{n}{m}]$, to the
$n$-skeleton of the Eilenberg-Maclane space $K(\Z,m)$ which induces an
epimorphism in rational homology in all dimensions through $n$.  We
construct a CW complex $W$ by attaching cells to $(S^m)^{\times p}$ so
as to replace the epimorphisms by isomorphisms.  By Sullivan \cite{S},
the localisation at 0 (rationalisation) $W_0$ of $W$ may be thought of
as the $n$-skeleton of $K(\Q,m)$.  Now a map from a compact complex
$X$ has image in a finite subcomplex of $K(\Q,m)$ and hence in a
finite piece of $W_0$.  If $W_0$ admits a telescope model, we can
conclude that the image of $X$ may be deformed into a copy of $W$
inside $W_0$.  In general, the argument is more delicate; here we need
to overcome two difficulties: 
\begin{alphenum}
\item
lack of control of the higher homotopy groups of the skeleta of
James's model of the spherical loop space (\cf Remark
\ref{difficult});
\item
the extra dimension of cells present in a `rational cell' being
attached in the process of rationalizing a complex (\cf Remark
\ref{difficult2}).
\end{alphenum}
Finally, the pullback arguments of \cite{B}, \cite{BKS}, and \cite{KS}
reduce the freedom of skeleta of the loop space to the freedom of
$(S^m)^{\times p}$, known since \cite{K}.

\begin{remark}[Pair of complementary dimensions]
\label{pair}
The $m$-systole $\sys_m(X)$ of $X$ can be briefly defined as the least
mass of a rectifiable $m$-current representing a nonzero class in
integer homology (\cf section 3).  Our main theorem states that $X$
admits a sequence of metrics $\{g_j\}$ with the following asymptotic
behavior as $j\to \infty$:
$${\vol_n(g_j) \over \sys_m(g_j)^{{n\over m}}}\to 0,
$$
at least if the $m$-dimensional homology of $X$ is torsion free.  

A related problem is that of systolic freedom in a pair of
complementary dimensions $m$ and $n-m$.  Recall that $X$ is called
$(m,n-m)$-systolically free if there exist metrics $g_j$ which behave
as follows:
$${\vol_n(g_j) \over \sys_m(g_j) \sys_{n-m}(g_j)}\to 0.
$$
While our technique relies on classifying spaces and rational homotopy
theory, the problem of a pair of complementary dimensions can be
solved by a direct geometric construction, with the input from
algebraic topology reduced to Poincar\'e duality.  See \cite{B2, B3,
K3}.

\end{remark}

The paper is organized as follows.  Section \ref{sec:i} discusses
systolic constraint and mentions a recent result of M. Freedman on
freedom with $\Z_2$-coefficients.  Section~\ref{sec:loop} defines
systolic freedom and the concept of an `$m$-meromorphic map'.
Sections \ref{sec:proof} and \ref{sec:proof2} contain the technical
statement of Theorem \ref{thm:m} and its proof, and also describe the
pertinent telescoping ideas from rational homotopy theory.
Section~\ref{sec:whitehead} contains some speculations regarding an
alternative proof using higher order Whitehead products.
Section~\ref{torsion} outlines a program of study of systolic freedom
with torsion coefficients, based on an analysis of `torsion
meromorphic maps'.

\section{Systolic freedom versus constraint}
\label{sec:i}
Our theorem shows that the total volume imposes no upper bounds on the
size of the least-area homologically nontrivial submanifold.  In
contrast, let us recall two results on the existence of such bounds
for related invariants.

The first result is M. Gromov's inequality \cite{G1} for essential
manifolds, which we will state in the particular case of real
projective space $\R P^n$.  Namely, all metrics $g$ on $\R P^n$ obey
the inequality
\begin{equation}
\label{eq:1}
\min_{\gamma\sim \R P^1} \hbox{$g$-\length}(\gamma) \leq C_n {\root n
\of {vol(g)}},
\end{equation}
where the minimum in the lefthand side is over all closed curves in
$\R P^n$ which are homotopic to $\R P^1\subset \R P^n.$ Here the
constant $C_n$ is independent of $g$, and the $n$-th root ensures
scale invariance.  Gromov's theorem generalizes classical results of
Loewner and Pu in dimension~$2$ (see also \cite{H2}).

The second result is an inequality similar to the above, but with
length replaced by area.  Namely, all metrics $g$ on complex
projective space $\C P^n$ obey Gromov's inequality (see
\cite[p.~262]{G2}, or \cite[Theorem 2.6]{KS}):
\begin{equation}
\label{eq:2}
\inf_{\gamma^{\phantom{a}}_\Q\sim \C P^1}
\hbox{$g$-area}(\gamma^{\phantom{a}}_\Q) \leq C_n {\root n \of
{vol(g)}},
\end{equation}
where the infimum is over all {\it rational} $2$-cycles
$\gamma^{\phantom{a}}_\Q=\sum_k r_k \sigma_k\ (r_k\in \Q)$
representing the generator $[\gamma^{\phantom{a}}_\Q]=[\C P^1]\in
H_2(\C P^n,\Q)=\Q$ (the sharp value of the constant is ${\root n \of
{n!}}$, as discussed in \cite[section 6.3]{M} and \cite[1.8.2]{Fe}).
Here by definition
\begin{equation}
\label{eq:25}
\hbox{$g$-area}({\textstyle\sum}_k r_k \sigma_k)= {\textstyle\sum}_k |
r_k | \hbox{$g$-area}(\sigma_k),
\end{equation}
and the area of a piecewise smooth singular $2$-simplex $\sigma_k$
is the integral of the pullback of $g$ to the standard $2$-simplex.
This inequality was generalized by J. Hebda \cite{H1}.

\begin{remark}[Torsion coefficients]
\label{tors}
All of our results ultimately rely on a calibration technique using
differential forms (\cf \cite[Appendix A]{KS}).  A new technique using
instead a lower bound on the first eigenvalue $\lambda_1$ of the
Laplacian, was recently invented by M. Freedman \cite{Fr}, who proved
that $S^1 \times S^2$ is $(1,2)$-systolically free even if one uses
$\ZZ$-coefficients in homology when defining the systoles.  This
result implies in particular that the $4$-manifold $S^2\times S^2$ is
free in middle dimension even with $\ZZ$-coefficients (\cf section 7).

On the other hand, the question whether inequality \eqref{eq:2}
remains valid for the complex projective plane if we work with
coefficients in $\Z_2$ is still open:
\begin{equation}
\label{eq:3}
\inf_{\gamma^{\phantom{a}}_\ZZ\sim \C P^1}
\hbox{$g$-area}(\gamma^{\phantom{a}}_\ZZ)
\displaystyle\mathop\leq^{??}
C_n {\root n \of {vol(g)}},
\end{equation}
where $\gamma^{\phantom{a}}_\ZZ$ is a $2$-cycle with
$\Z_2$-coefficients (\eg a possibly nonorientable surface)
representing the nonzero element in $H_2(\C P^n,\ZZ)=\ZZ$.  
See also section \ref{torsion}.

It is still unknown whether $\R P^3$ is (1,2)-systolically free modulo
2 (this question was originally posed \cite[p.\ 622]{BBK} in 1994), or
in formulas:
\begin{equation}
\label{eq:}
\inf_{ \gamma^{\phantom{a}}_\ZZ,\: \sigma^{\phantom{a}}_\ZZ }
\hbox{$g$-length}(\gamma^{\phantom{a}}_\ZZ)
\hbox{$g$-area}(\sigma^{\phantom{a}}_\ZZ)
\displaystyle\mathop\leq^{??}  \text{Const}\: {vol(g)},
\end{equation}
for all metrics $(\R P^3,g)$, where the infimum is over all curves
$\gamma^{\phantom{a}}_\ZZ\sim \R P^1$ and surfaces $\sigma
^{\phantom{a}} _\ZZ\sim \R P^2$.
\end{remark}
\section{Pulling back metrics by $m$-meromorphic maps}
\label{sec:loop}
In \cite{KS}, we reduced the proof of the $m$-systolic freedom modulo
torsion (see definition below) of all smooth $n$-manifolds, to the
examination of a finite list of objects: the successive skeleta (up to
the $n^{th}$ skeleton) of the based loop space $\Omega(S^{m+1})$.  The
price one has to pay is the enlargement of the category to that of
CW complexes (see Definition \ref{def:systfree}).

\begin{definition}
\label{def:v}
Let $(X,g)$ be a finite $n$-dimensional simplicial complex,
endowed with a piecewise smooth metric $g$.  Let
$m\leq n$.  Let $\alpha\in H_m(X,\Z)$.  Define
\[
v(\alpha)=\inf_{M\in\alpha}\vol_m(M),
\]
where the infimum is taken over all piecewise smooth integer cycles
$M$ representing the class $\alpha$.  Here the volume of a (smooth)
singular simplex is that of the pullback of the quadratic form $g$ to
the simplex, and we take absolute values of the coefficients to obtain
the volume of the cycle, as in formula \eqref{eq:25}.
\end{definition}

We define the {\em systole modulo torsion}, $\sys^\infty$, of $(X,g)$
by setting
\[
\sys^{\infty}_m(g)=\inf_{\substack{\alpha\in H_m(X,\Z)\\
\alpha\ne{\rm torsion}}}v(\alpha).
\]

\begin{definition}
\label{def:systfree}
An $n$-dimensional CW complex $X$ is {\em
$m$-systolically free (modulo torsion)} if
\[
\inf_{g} \frac{\vol_{n}(g)^{\frac{m}{n}} }{ \sys^{\infty}_m (g)}=0,
\]
where the infimum is taken over all piecewise smooth metrics $g$ on a
finite simplicial complex $X'$ homotopy equivalent to $X$.
\end{definition}

Here the choice of $X'$ is immaterial by virtue of the pullback Lemma
\ref{lem:xyfree} below, which applies in particular to homotopy
equivalences.  The idea of the proof of the systolic freedom of the
skeleta of the loop space is a reduction to the case of the product of
spheres, for which systolic freedom was established in \cite{K} for
$m\geq 3$, and in \cite[Lemma 4.5 and Corolary 7.9]{KS} in the
remaining case:

\begin{prop}[\cite{K}, \cite{KS}]
\label{prop:spheres}
The Cartesian product $(S^m)^{\times \l}$ of $\l$ copies of the
$m$-sphere is $m$-systolically free, for all $m\geq 2$ and $\l\geq 2$.
\end{prop}

We will carry out such a reduction by means of constructing a morphism
from a skeleton of $\Omega(S^{m+1})$ to a product of spheres, and
applying pullback arguments developed in \cite{B}, \cite{BKS}, and
\cite{KS}.

The morphisms from $X$ to $Y$ best suited to our problem are more
flexible than continuous maps from $X$ to $Y$.  We still work with
continuous maps, but we allow certain enlargements of the target $Y$
of the map, which we will refer to as `meromorphic extensions'.

\begin{definition}
\label{def:meroe}
Let $Y$ be an $n$-dimensional CW complex.  A CW complex $W$ is called
a {\em meromorphic extension} of $Y$ if $W$ has the homotopy type of a
finite CW complex obtained from $Y$ by attaching cells $e^d$ of
strictly smaller than the top dimension, or in formulas:
$$W\simeq Y\cup \bigcup_{i} e^{d_i},\text{ where } d_i\leq n-1.$$
\end{definition}

\begin{definition}
\label{def:mero}
Let $X^n$ and $Y^n$ be finite CW complexes of dimension $n$.  An {\em
$m$-meromorphic map}, $X\mero Y$, is a pair $(W,f)$ where $W$ is a
meromorphic extension of $Y$ and $f:X\to W$ is a continuous map which
induces a monomorphism in $m$-dimensional rational homology:
$$f_*: H_m(X,\Q) \to H_m(W,\Q)\text{ is injective.}$$
\end{definition}

Our terminology is inspired by the observation that the blow-up map
$Y=\hat X\rightarrow X$ of, say, a complex analytic manifold $X$
admits a kind of an inverse: $X\mero Y$.  Here the space $W$ is
obtained by coning off the exceptional divisor in $Y=\hat X$, while
the inverse is perturbed in a neighborhood of the blown up point (see
\cite[Example 7.4]{KS} for details).  Such maps are useful because of
the following lemma on pulling back systolic freedom (\cf
\cite[Proposition 7.10]{KS}).

\begin{lemma}[Pullback Lemma]
\label{lem:xyfree}
Let $X \mero Y$ be an $m$-meromorphic map.  If $Y$ is $m$-systolically
free modulo torsion, then so is $X$.
\end{lemma}

\begin{example}
\label{ex}
Let us show how one can apply the pullback lemma to reduce the
2-systolic freedom of the complex projective plane $\C P^2$ to that of
the product $S^2\times S^2$.

Consider the decomposition $\C P^2=S^2 \cup_h D^4$, where
$h={\frac{1}{2}[e,e]}$ is the Hopf map.  Here $e$ is the identity map
of $S^2$, and $[e,e]$ is the Whitehead product.  Also consider the
CW complex $W=S^2\times S^2\cup_{e_1 -e_2}D^3$.  Let $f:\C P^1\to W$,
$f(e)=de_1$ where $d\geq 2$ is an even integer. Then
\begin{equation*}
\begin{split}
f(h)=f(\tfrac{1}{2}[e,e])=\tfrac{1}{2}[f(e),f(e)]= \tfrac{1}{2}
[de_1,de_1] \\ =\tfrac{d^2}{2}[e_1,e_1] =\tfrac{d^2}{2}[e_1,e_2]=0
\text{\ in\ } \pi_3(W),
\end{split}
\end{equation*}
since $e_1=e_2$ in $W$, and so $[e_1,e_1]= [e_1,e_2]$.  Thus $f$
extends to a map $f: \C P^2\to W$ which is injective in 2-dimensional
homology, and hence may be viewed as a 2-meromorphic map $\C P^2\mero
S^2 \times S^2$.  Now Lemma \ref{lem:xyfree} reduces the 2-systolic
freedom of $\C P^2$ to that of $S^2\times S^2$.

See also section~\ref{sec:whitehead} for a possible generalisation to
$\C P^n$ using higher order Whitehead products.
\end{example}

The above argument is valid in view of the absence of torsion in
$\pi_{2m-1}(S^m)$ for $m=2$.  For $m>2$, we use the following argument
to find an $m$-meromorphic map to $S^m\times S^m$ (see proof of Lemma
4.4 in \cite{KS}).

\begin{example}
\label{ex2}
Let $X=S^m\cup_a D^{2m}$ be a CW complex of dimension $2m$ with 3
cells of dimensions 0, $m$, and $2m$, respectively.  Let $W=S^m\times
S^m\cup_{e_1-e_2}D^{m+1}$.  Recall that the Whitehead product
$[e,e]\in \pi_{2m-1}(S^m)$ generates precisely the kernel of the
suspension homomorphism.  Suspension commutes with the homomorphism
induced by the degree $q$ self-map $\phi_q: S^m\to S^m$.  Hence, if
$q$ is a multiple of the order of the (finite!)  stable group
$\pi_{2m}(S^{m+1})$, then the map $\phi_q$ sends $\pi_{2m-1}(S^m)$ to
the subgroup generated by Whitehead products.  In particular,
$\phi_q(a)=\lambda [e,e]$.  We now compose $\phi_q$ with the inclusion
$\iota:S^m\to W$ of $S^m$ as the first factor of $S^m\times S^m\subset
W$.  It follows that
$$\iota\circ\phi_q(a)=\iota(\lambda [e,e])=\lambda [e_1,e_1]=
\lambda[e_1,e_2]= 0\in \pi_{2m-1}(W).$$ 
Thus, the map $\iota\circ\phi_q$
extends to a map $X\to W$.  This proves the $m$-systolic freedom of
$X$.
\end{example}

\section{Proof of systolic freedom of loop space}
\label{sec:proof}
The results of this section generalize the middle-dimensional systolic
freedom, established in \cite{BK}, \cite{BKS}, and \cite[Theorem
1.5]{KS}.  Let $n=mp,$ where $m\geq 2$ and $p\geq 2$.  Let
$(S^m)^{\times p}$ denote the $p$-fold Cartesian product of
$m$-spheres.

\begin{theorem}
\label{thm:main}
Let $X$ be a finite $n$-dimensional CW complex with unit $m$-th Betti
number: $b_m(X)=1$.  Then $X$ admits an $m$-meromorphic map to
$(S^m)^{\times p}$, and hence $X$ is $m$-systolically free modulo
torsion.
\end{theorem}

\begin{proof}
Let $n = mp$, where $m$ is even (the case of $m$ odd is easier; see
\cite[section~10]{KS}).  Let $S=(S^m)^{\times p}$.  Let $\phi_j$ be
the selfmap of $S$ defined by a degree $j$ map on each factor.  The
map $\phi_j$ induces a {\em scalar} homomorphism (namely,
multiplication by a power of $j$) in each cohomology group of $S$.
Namely, $\phi_j^*=\wedge^\ell(j\id)=j^\ell: H^{\ell m} (S) \to H^{\ell
m}(S).$ Let $w = x_1 +\dots + x_p$ be the sum of the standard
generators $x_i$ of the cohomology group $H^m(S)=\Z^p$.  Then
$$w^p=p!\: x_1\cup\cdots\cup x_p \ne 0 \text{\ in\ }H^{n}(S).$$ Now
the class $w$ defines a map $f: S \rightarrow K(\Z,m)$ such that
$f^*(u) = w,$ where $u\in H^m(K(\Z,m))=\Z$ is a generator.  Recall
that $u$ is a multiplicative generator of the polynomial ring
$H^*(K(\Z,m)) \otimes \Q\cong \Q[u]$ ($m$ even).  Our proof will rely
on the following lemma.

\begin{lemma}
\label{lem1}
The map $f: (S^m) ^{\times p} \rightarrow K(\Z,m)$ extends to a map
$f:W\to K(\Z,m)$ which induces an isomorphism in rational homology
through dimension $n$, where $W$ is a finite CW complex obtained from
$S$ by attaching cells of dimension at most $n-m+1<n$.  Furthermore,
the maps $\phi^2_j$ extend to $W$ for $j$ divisible by a sufficiently
large $j_0$, in the following two situations:
\begin{alphenum}
\item if $p=2$ or $3$;
\item if $m=2$.
\end{alphenum}
\end{lemma}

\begin{proof}
Since $w^p \not = 0$, the map $f: (S^m) ^{\times p} \rightarrow
K(\Z,m)$ induces an epimorphism
\begin{equation}
\label{epi}
H_{k}(S,\Q)\rightarrow H_{k}(K(\Z,m),\Q)
\end{equation}
in each dimension less than $n$, and an isomorphism in dimension $n$.
We will use the relative Hurewicz theorem to eliminate the successive
kernels of the homomorphisms \eqref{epi}, and construct $W$ by skeleta
$W^{(k)}$ inductively on $k$ (\cf \cite[p.~27]{S}).

Let $S=(S^m)^{\times p}$ and let $W^{(m)}=S^{(m)}=S^m\vee\cdots\vee
S^m$.

Next, let $W^{(m+1)}=S^{(m)}\cup_{e_1-e_2}D^{m+1} \cup_{e_2-e_3} \cdots
\cup_{e_{p-1}-e_p} D^{m+1}$.  Then $W^{(m+1)}\simeq S^m$.

Since the map $\phi_j$ acts by $j\id$ on $\pi_m(W^{(m)})= H_m(S)$, it
extends to $W^{(m+1)}$.

We set $W^{(2m)}= W^{(m+1)}\cup S^{(2m)}$ and note that, by patching
the two pieces, $\phi_j$ extends to $W^{(2m)}$.  Note that
$W^{(2m)}\simeq S^m\cup_i D^{2m}_i$, where $i=1,\ldots,
\tbinom{p}{2}$.  We have $W^{(2m)}/W^{(m+1)}\simeq \vee_i S^{2m}$,
$i=1,\ldots, \tbinom{p}{2}$.  Consider the exact sequence
\begin{equation}
\label{seq}
\pi_{2m}(W^{(m+1)})\to
\pi_{2m} (W^{(2m)})\to \pi_{2m}(W^{(2m)}, W^{(m+1)}).
\end{equation}
By Hurewicz's theorem, the map $\phi_j$ induces multiplication by
$j^2\id$ on the group
\begin{equation}
\label{cannot}
\pi_{2m}(W^{(2m)}, W^{(m+1)})= H_{2m}(W^{(2m)}, W^{(m+1)})= H_{2m}
(\vee_i S^{2m})=\Z^{\tfrac{p}{2}}_{\phantom{a}}.
\end{equation} 
Now consider the diagram
\begin{equation*}
\label{eq:cd1}
\begin{CD}
H_{2m+1}(K(\Z,m),W^{(2m)}) @>\partial>>
H_{2m}(W^{(2m)}) @>f_*>> H_{2m}(K(\Z,m))  \\ 
@AA{h_1}A
@AA{h_2}A \\ 
\pi_{2m+1}(K(\Z,m),W^{(2m)}) @>\partial>>
\pi_{2m}(W^{(2m)}),
\end{CD}
\end{equation*}
where $h_i$ are the Hurewicz maps.  By the relative Hurewicz theorem,
the homomorphism $h_1$ is surjective.  Thus
$\ker(f_*)=\im(\partial\circ h_1) =\im(h_2\circ\partial).$ Hence,
every element in $\ker(f_*)$ is spherical.

Let $\alpha_i\in \pi_{2m}(W^{(2m)})$, $i=1,\ldots \tfrac{p}{2}-1$, be
a set of lifts via the Hurewicz homomorphism $h_2$, of a set of
generators of the kernel of $f$.  From the exact sequence \eqref{seq},
we have $\phi_j(\alpha_i)= j^2 \alpha_i+t$, where $t\in\pi_{2m}(S^m)$
is of finite order.  By Sullivan's result \cite[p.\ 19]{S} on the
nilpotence of the finite homotopy groups of spheres, there exists a
$j_0$ such that if $j_0 | j$ then $\phi_j=0$ on $\pi_{2m}(S^m)$.
Assume also that $j_0$ is a multiple of the order of $\pi_{2m}(S^m)$.
Then
$$\phi_j^2(\alpha_i)=\phi_j(j^2\alpha_i+t)= j^2(\phi_j(\alpha_i))+
\phi_j(t)= j^4 \alpha_i+ j^2 t= j^4\alpha_i,$$ \ie $\phi_j^2
(\alpha_i)$ is proportional to $\alpha_i$.

It follows that $\phi_j^2$ extends to the next skeleton
$$W^{(2 m+1)}= W^{(2 m)}\cup_{\alpha_i} D_i^{2 m+1}
$$ 
and therefore also to $W^{(3m)}= W^{(2m+1)}\cup S^{(3m)}$, proving
part (a) of the lemma.

We can similarly define $W^{(\ell m)}=W^{((\ell-1)m+1)}\cup S^{(\ell
m)}$, and use the Serre form of the relative Hurewicz theorem (see
\cite[pp.~95--98]{MT}) to define
$$W^{(\ell m+1)}= W^{(\ell m)}\cup_{\alpha_i} D_i^{\ell m+1} \text{
for } i=1,\ldots,\tbinom{p}{\ell}-1.
$$ 
Proceeding by induction, we extend $f$ to the space $W=W^{((p-1)m+1)}
\cup S$, which is the desired meromorphic extension of $S$ (\cf
Definition \ref{def:meroe}).
\begin{remark}
\label{difficult}
The above argument does not yield an extension of $\phi_j^2$ to
$W^{(3m+1)}$ (when $p\geq 4$) because already the group $\pi_{3 m}
(W^{(2m+1)}, W^{(m+1)})$ cannot be shown to be finite merely from
Hurewicz's theorem, as in the calculation \eqref{cannot} above.  An
alternative approach given below gives an easy proof only for $m=2$
(\cf Remark \ref{hard}).
\end{remark}
For $m=2$, we can understand the effect of $\phi_j$ on $\pi_{2\ell}
(W^{(2\ell)})$ easily using the Hopf fibration over $\C P^n$.  We
argue inductively.  To the extent that $\phi_j$ already acts on
$W^{(2\ell)}$, the family $(\phi_j)$ localizes homology of the space
$W^{(2l)}$.  Therefore by \cite[p.\ 19]{S}, the family also localizes
homotopy, so that we have the same corollary as for the spheres: the
map on $d$ torsion of $\pi_i(W^{(2l)})$ induced by $\phi_d$ is
nilpotent.  To show that the action of $\phi_j$ is scalar on
$\pi_{2\ell} (W^{(2\ell)})$, it therefore suffices to show that this
group has rank precisely
$$\text{rank}(H_{2\ell} (W^{(2\ell)}))-1=\tbinom{p}{\ell}-1.$$
Consider the exact sequence of pair $(W^{(2\ell)},W^{(2l-1)})$:
$$\pi_{2\ell} (W^{(2\ell -1)}) \to \pi_{2 \ell}(W^{(2\ell)})\to \pi_{2
\ell} (W^{(2\ell)}, W^{(2\ell-1)})= H_{2 \ell} (W^{(2\ell)},
W^{(2\ell-1)})= \Z^{\binom{p}{\ell}}$$ by Hurewicz.  We need to show
that the group $\pi_{2\ell} (W^{(2\ell -1)})$ is finite.  Now the map
$W^{(2\ell -1)}\to \C P^{\ell-1}$, where $\C P^{\ell-1} \subset \C
P^\infty=K(\Z,2)$, has already been constructed.  This map induces an
isomorphism in rational homology in all dimensions.  By the
Serre-Hurewicz theorem, the group $\pi_i(\C P^{\ell -1}, W^{(2\ell
-1)})$ is finite for all $i$.  The exact sequence of the pair shows
that the groups $\pi_{i}(W^{(2\ell-1)})$ are all finite except $\pi_2$
and $\pi_{2\ell -1}$, since $\pi_{2\ell-1}(\C
P^{\ell-1})=\pi_{2\ell-1}(S^{2\ell-1})=\Z$ from the exact sequence of
the Hopf fibration, proving part (b) of the Lemma.
\end{proof}
\begin{remark}
\label{hard}
To generalize this proof to $m\geq 3$, one could replace the spaces
$\C P^n$ by the skeleta of $\Omega(S^{m+1})$ (\cf formula
\eqref{james}).  However, the calculation of the homotopy groups of
such skeleta is not as immediate as that of complex projective spaces.
\end{remark}
\section{Conclusion of proof}
\label{sec:proof2}
Lemma \ref{lem1} provides a quick proof of Theorem \ref{thm:main} in
the cases (a) and (b), as follows.  We appeal to rational homotopy
theory, to conclude that the rationalisation $W_0$ of $W$ can be
thought of as the $n$-skeleton of $K(\Q,m)= K(\Z,m)_0$.

Due to the existence of the selfmaps $\phi^2_j$ of $W$, according to
D. Sullivan \cite{S}, the space $W_0$ admits a model as an infinite
telescope on $W$, whose $j$-th stage is the cylinder of the self-map
$\phi^2_j$ of $W$, where $\phi_j$ is induced by a degree $j$ map on
each of the factors in $S^m\times \cdots \times S^m=S$.

Now let $X$ be a finite $n$-dimensional CW complex with
$b_m(X)=1$. Consider a map from $X$ to $K(\Q,m)$ defined by any
non-torsion $m$-dimensional cohomology class.  The map may be assumed
to have image in $W_0$, viewed as the skeleton of the classifying
space.  Being compact, the image of $X$ lies in a finite piece of the
telescope structure of $W_0$.  Hence it can be retracted to the final
stage, $W$, of the finite piece.  Now $W$ is just the space
$(S^m)^{\times p}$, with some cells of dimension at most $n-m+1<n$
attached.  Hence $W$ is a meromorphic extension of $S$.  We thus
obtain an $m$-meromorphic map $X\mero (S^m)^{\times p}$.  The pullback
lemma \ref{lem:xyfree} completes the proof in the 2 cases mentioned in
Lemma \ref{lem1}.

In the general case $p\geq 4, m\geq 3$, we need a more delicate
argument.  We rely on the following property of maps of compact spaces
into rationalisations.

\begin{lemma}
\label{lem9}
Let $S$ be a CW complex admitting a telescope model (\ie there are
self-maps which localize homology).  Assume that a CW complex $W$ is
obtained from $S$ by attaching cells of dimension at most $k$.  Then
the image of a compact space mapping into the rationalisation $W_0$ of
$W$ may always be deformed into a subcomplex of $W_0$ which is
homotopy equivalent to a copy of $S$ with cells of dimension at most
$k+1$ attached.
\end{lemma}

\begin{proof}
Recall that $S_0$ is obtained as the direct limit in the following
construction:
\begin{equation}
\label{frombott}
S_0=\displaystyle\mathop\amalg_j S \times I\ /\ (x,1)\sim (\phi_j
\circ id_j(x),0).
\end{equation}
Here $id_j: S\rightarrow S$ is the identification
of the $j$-th and $(j+1)$-th levels.

A local CW complex is built inductively by attaching cones over the
local sphere using maps of the local spheres into the lower `local
skeletons'.  For each cell of dimension $\leq k$ attached to $S$, we
attach a corresponding local cell to $S_0$, which contains cells
$\sigma^{d_i}$ of dimension $d_i\leq k+1$.  Here the extra dimension
is due to the presence of cylinders in formula \eqref{frombott}.  A
map from a compact space $X$ into $W_0$ has image in a finite
subcomplex $W'_0\subset W_0$, which may be assumed to be of the
following form: take a finite piece
$$S'_0=\displaystyle\mathop\amalg_j^N S \times I\ /\ (x,1)\sim (\phi_j
\circ id_j(x),0)$$ of the infinite telescope $S_0$, and attach, to
$S_0'$, at most finitely many cells from among the $\sigma^{d_i}$.
Namely,
$$W'_0= S'_0\cup_i \sigma^{d_i}.$$ From the homotopy equivalence
$S'_0\simeq S$, we obtain the equivalence $W'_0\simeq S\cup_i
\sigma^{d_i}$, since each attaching map is an element of a homotopy
group and hence a homotopy invariant, while homotopy groups commute
with direct limits by a simple compactness argument.
\end{proof}
To complete the proof of Theorem \ref{thm:main}, we argue as follows.
Let $W_0$ be the rationalisation of the space $W$ of the lemma.  By
the universal property of a localisation \cite[p.\ 18]{S}, the map
from $W$ to the local space $K(\Q,m)$ induced by $f$ extends to a map
$W_0\to K(\Q,m)$.  Thus $W_0$ may be thought of as the $n$-skeleton of
$K(\Q,m)$.  Recall that $W$ is obtained from $S$ by adding cells of
dimension at most $n-m+1$.

If $X$ is compact, by Lemma \ref{lem9} the image of the map can be
deformed into a space $W'_0$ obtained from $S$ by attaching cells of
dimension at most $n-m+2$.  
\begin{remark}
\label{difficult2}
The space $W'_0$ is thus a meromorphic extension (\cf Definition
\ref{def:meroe}) of $S$ only if $m\geq 3$.  The case $m=2$ is handled
as above (following Remark \ref{difficult}), using the Hopf fibration
over $K(\Z,2)=\C P^\infty$.
\end{remark}
Since we rely on an existence theorem for rationalisations, what we
lose control of in this version of the proof is the exact form of the
meromorphic extension of $S$, which admits a continuous map from $X$.
\end{proof}

\begin{cor} 
\label{coro}
Every finite $n$-dimensional CW complex is $m$-systolically free
modulo torsion, provided $2\leq m <n$.

\end{cor}

\begin{proof} 
Recall that Theorem 1.6 of \cite{KS} reduces the general case to that
of the successive skeleta of the loop space $\Omega=\Omega(S^{m+1})$,
arising from the cell decomposition of I. James,
\begin{equation}
\label{james}
\Omega \simeq S^m \cup e^{2m} \cup e^{3m} \cup \cdots ,
\end{equation}
with precisely one cell in each dimension divisible by $m$ (see
\cite{Wh}).  We now apply Theorem \ref{thm:main} to the $n$-skeleton
of $\Omega$, to prove its $m$-systolic freedom modulo torsion,
completing the proof of the Corollary.

Our proof of Theorem 1.6 of \cite{KS} can be simplified by using the
telescope model for the localisation $\Omega_0=K(\Q,m)$ of the loop
space $\Omega= \Omega (S^{m+1})$.  Such a model exists by \cite{S},
since $\Omega$ admits self-maps defined by mapping a loop to its
$j$-th iterate, which induce multiplication by $j$ in all homology
groups.

Namely, let $X$ be a finite $n$-dimensional CW complex, and let
$b=b_m(X)$ be its $m$-th Betti number.  A choice of a basis for a
maximal lattice in $H^m(X)$ defines a map $X\to K(\Z^b,m)$.  The
composition of this map with the canonical map $K(\Z^b,m)\to
K(\Z^b,m)_0= \Omega_0^{\times b}$, induces an isomorphism $H_m(X,\Q)
\to H_m(\Omega_0^{\times b},\Q)$.  Now we argue as above.  We use the
compactness of $X$ to construct a projection to the final level of a
finite piece of the telescope.

By the pullback lemma \ref{lem:xyfree}, it remains to prove the
$m$-systolic freedom of the $n$-skeleton of $\Omega^{\times b}$.
Next, we apply the carving up technique of \cite{KS}, to reduce the
problem to the $m$-systolic freedom of the closures of the
top-dimensional cells in the $n$-skeleton $(\Omega^{\times b})^{(n)}$
of $\Omega^{\times b}$.

The cell decomposition of $(\Omega^{\times b})^{(n)}$ induced by the
James decomposition \eqref{james} contains only cells of dimensions
divisible by $m$.  In particular, there are no cells of codimension~1.

The absence of codimension 1 cells in $(\Omega^{\times b})^{(n)}$ is
the crucial ingredient of the carving up technique.  It allows us to
isolate the different cell closures in $(\Omega^{\times b})^{(n)}$
from each other.  This is accomplished by inserting long cylinders
based on the {\em boundary}, in $(\Omega^{\times b})^{(n-m)}$, of each
top-dimensional cell.  The absence of cells of codimension 1
guarantees that these cylinders have positive codimension, and thus
make no contribution to total $n$-dimensional volume.

The effect of the long cylinders is to minimize the interaction
between distinct top-dimensional cells, so that an $m$-cycle traveling
from one to the other would have to pay a heavy price in terms of its
$m$-volume.  The formal argument, using coarea and isoperimetric
inequalities, appears in~\cite[Appendix B]{KS}.

Each of the top-dimensional cells in $(\Omega^{\times b})^{(n)}$ is a
product of skeleta $\Omega^{(mp_i)}$ of $\Omega$.  Each
$\Omega^{(mp_i)}$ admits a meromorphic map to $(S^m)^{\times p_i}$ by
Theorem \ref{thm:main}.  Thus there exists a meromorphic map
$$\textstyle{\prod}_i\:\Omega^{(mp_i)}\mero\textstyle{\prod}_i\:
(S^m)^{\times p_i}= (S^m)^{\times({\Sigma}_i p_i)} =(S^m)^{\times p}.
$$
The $m$-systolic freedom of $\textstyle{\prod}_i\:\Omega^{(mp_i)}$ now
follows from the systolic freedom of products of spheres (\cf
Proposition~\ref{prop:spheres}).
\end{proof}

\section{Higher order Whitehead products}
\label{sec:whitehead}
It is tempting to try to generalize the argument of Example \ref{ex}
using higher order Whitehead products, so as to derive the systolic
freedom of $\C P^n$, or of the appropriate skeleton of the loop space,
from the systolic freedom of the Cartesian product $(S^2)^{\times n}$.

We form the analogue of the space $W$ of Example \ref{ex} by defining
$W_n$ to be the CW complex obtained from $(S^2)^{\times n}$ by
attaching $n-1$ copies of a 3-cell along differences of consecutive
2-spheres in the product $(S^2)^{\times n}$:
$$W_n=(S^2)^{\times n}\cup_{e_1-e_2}D^3\cup_{e_2-e_3} \cdots
\cup_{e_{n-1}-e_n} D^3.$$ The top-dimensional cell of $(S^2)^{\times
n}$ is attached by means of the $n^{\text{th}}$ order Whitehead
product $a_n=[e_1,\ldots,e_n]$ of the 2-dimensional generators.  Thus
$$W_{n}=W_n^{(2n-2)}\cup_{a_n} D^{2n}$$ where $W^{(i)}$ denotes the
$i$-skeleton.  Furthermore,
$\C P^n=\C P^{n-1}\cup_h D^{2n}$ where 
\begin{equation}
n!\:h=[e,e,\dots ,e]
\label{eq:hain}
\end{equation}
(\cf \cite[p.~416]{P}, \cite[p.~460]{A}, \cite[p.~39]{Ha}).  Now let
$f=\phi_d$ be a self-map of $S^2$ of degree $d$.  By virtue of formula
\eqref{eq:hain}, we can write
$$
f(n!\:h)=f([e,e,\ldots,e])\in [f(e), f(e),\ldots, f(e)]=
[de_1,de_1,\ldots,de_1].
$$
Since $e_i=e_j$ in $W_n^{(2n-2)}$ for all $i$ and $j$, we have
$$f(n!\:h) \in [de_1,de_2,\ldots,de_{n}] = d^n [e_1,e_2,\ldots,e_n].
$$ Since $f(n!\:h)= n!\:f(h)$, we obtain
\[
f(h)\in\frac{d^n}{n!} [e_1,e_2,\ldots,e_n]+ t\subset \pi_{2n-1}
(W_n^{(2n-2)})
\]
assuming the fraction is an integer (\eg if $n!$ divides $d$), where
$t\in \pi_{2n-1}(W_n^{(2n-2)})$ is a torsion element of order dividing
$(n!)$.

Now the problem of the torsion element $t$ can be handled by replacing
$\C P^n$ by the $2n$-skeleton of a model of $\Omega(S^3)$ (\cf
\eqref{james}).  Then the attaching map $h'$ is a higher order
Whitehead product and we may set $d=1$, to obtain
$$h'\in [e_1,\ldots,e_n] \subset\pi_{2n-1}(W_n^{(2n-2)}).
$$ Now the set $[e_1,\ldots,e_n]\subset\pi_{2n-1}(W_n)$ clearly
contains zero, corresponding to the attaching map $a_n$ of the
top-dimensional cell in the product of spheres.  However, this class
may also contain other elements, due to ambiguity in the choice of the
extension to the fat wedge.  For instance, the 4-skeleton of the space
$W_3=S^2\times S^2\times S^2 \cup_{e_1-e_2}D^3 \cup_{e_2-e_3}D^3$ is
homotopy equivalent to the 2-sphere with three distinct 4-cells
attached along the same class $[e,e]$.  There are thus 9 potentially
nonhomotopic possibilities for extending the inclusion $S^2\vee
S^2\vee S^2\to W_3$ to the fat wedge $W_3^{(4)}\to W_3$.

Thus, the indeterminacy of the higher order Whitehead product is an
obstruction to the attempted generalisation of the proof of Example
\ref{ex}.  On the other hand, attaching higher-dimensional cells to
reduce the rank of all homology groups to at most 1 may resolve the
indeterminacy and yield an alternative proof of the main theorem.

\section{Torsion coefficients for meromorphic maps}
\label{torsion}

There now exists a first example of systolic freedom even with $\ZZ$
coefficients, due to M. Freedman \cite{Fr} (\cf Remark \ref{tors}).
One could therefore propose to study systolic freedom over $\ZZ$ of
manifolds.  Essentially this would amount to proving the statement of
Theorem \ref{thm:m}, with the word `orientable' deleted, as in the
following question.
\begin{question}
\label{quest}
Let $m$ and $n$ be integers such that $2\leq m<n$.  Let $X$ be an
$n$-dimensional compact smooth manifold.  Does $X$ admit metrics of
arbitrarily small total volume, such that every $m$-dimensional
(possibly nonorientable) submanifold of less than unit $m$-volume is
null-homologous as a cycle with $\ZZ$ coefficients?
\end{question}
We could attack this question by studying `torsion meromorphic maps',
where we replace $\Q$ by $\ZZ$ in Definition \ref{def:mero}.  We could
then study the partial order on the set of all manifolds of a given
dimension, defined by the existence of a torsion meromorphic map
between them.

This approach immediately yields the 2-systolic freedom over $\ZZ$ of
$S^2 \times S^2$, since it admits a degree 1 map to $S^2\times
S^1\times S^1$ with a pair of 2-cells attached (\cf \cite[Corollary
7.9]{KS}), while the freedom of $S^2\times S^1\times S^1$ follows from
Freedman's example $(S^2\times S^1,g_F)$ by forming Cartesian product
with a circle of length $\frac{\sys^{\ZZ}_2(g_F)}{\sys^{\ZZ}_1(g_F)}$.

Meanwhile, a simple analysis of the multiplicative structure of the
pertinent cohomology rings suggests that the manifold $\C P^2$ ought
not to admit a torsion meromorphic map to $S^2\times S^2$.  It is thus
the simplest manifold for which Question \ref{quest} is open.

\begin{ack}
We are grateful to Shmuel Weinberger and to Emanuel Dror-Farjoun for
helpful discussions of rational homotopy theory.  Ivan Babenko
carefully went over an early version of this paper and made a number
of pertinent suggestions and corrections, for which we thank him.
\end{ack}

\bibliographystyle{amsalpha}

\end{document}